\documentclass[12pt]{amsart}
\usepackage{amsmath,amsfonts,euscript,amscd,amsthm,amssymb,upref,graphics}
\usepackage[retainorgcmds]{IEEEtrantools}
\theoremstyle{plain}
\newtheorem{theorem}{Theorem}
\swapnumbers

\newtheorem{proposition}[subsection]{Proposition}
\newtheorem{lemma}[subsection]{Lemma}
\newtheorem{corollary}[subsection]{Corollary}

\newtheorem{convention}[subsection]{Convention}

\theoremstyle{definition}
\newtheorem{definition}[subsection]{Definition}
\newtheorem{example}[subsection]{Example}
\newtheorem{remark}[subsection]{Remark}

\newtheorem{nothing*}[subsection]{}

\newcommand{\rien}[1]{}

\newcommand{\Spec}{ \operatorname{{\rm Spec}}}

\newcommand{\Aut}{ \operatorname{{\rm Aut}}}

\newcommand{\C}{\ensuremath{\mathbb{C}}}

\newcommand{\Ker}{{\rm Ker} \,}

\newcommand{\nn}{\nonumber}

\def\aut{{\rm Aut}}
\def\saut{{\rm SAut}}

\renewcommand{\epsilon}{\varepsilon}
\renewcommand{\phi}{\varphi}

\begin{document}
\title{Makar-Limanov invariants, Derksen invariants, flexible points}

\author{Fabrizio Donzelli}
\address{Institute of Mathematical Sciences, Stony Brook University, 
Stony Brook, NY, 11794 \ \ USA}
\email{fabrizio@math.sunysb.edu}

{\renewcommand{\thefootnote}{} }

\begin{abstract}
We estabilish some connections among the Makar-Limanov invariant, the Derksen invariant, and the existence of flexible points on an affine variety.
\end{abstract}

\maketitle \vfuzz=2pt

\section{Flexibility}

In \cite{AFKKZ} the authors explore the properties of the automorphism group affine varieties having flexible points, which are defined as follows.

\begin{definition}\label{flexible}
Let $X$ be an algebraic variety. A point $x\in X$ is said to be flexible 
if the span of the vectors at $x$ tangent to the orbits of $\C_+$-actions coincides with $T_xX$. $X$ is said to be flexible if all points of $X_{reg}$, the regular locus of $X$, are flexible.
\end{definition}

Let $\aut (X)$ denote the group of automorpshism of an affine variety $X$. We are interested in the subgroup
of $\aut (X)$ generated by $\C_+$-actions.

\begin{definition}\label{special}
For an algebraic variety $X$, the special automorphism group of $X$ is the subgroup of $\aut (X)$ generated by all subgroups of $\aut (X)$
isomorphic to $\C_+$ (that is all algebraic $\C_+$-actions). It will be denoted by $\saut (X) $.
\end{definition}

Since the flexibility of a point is an open condition, one expects
that if $X$ is flexible, then $\saut (X)$ acts transitively  on the regular locus of $X$, but a much stronger statement holds.

\begin{theorem}\label{ft}( \cite{AFKKZ}, Theorem 0.1)
Let $X$ be a flexible affine variety of dimension at least $2$. 
Then $\saut (X)$  acts $k$-transitively on the regular locus of $X$, for any integer $k>0$:
given two collections $\{x_1,\dots, x_k\}$ and $\{y_1,\dots y_k\}$
of pairwise dinstinct points, there exists a special automorphism $\psi\in \saut (X)$ such that  $\psi (x_i)=y_i$, for all $i\in\{1,\dots ,k \}$.
\end{theorem}

In the next sections we will state and apply more theorems from 
\cite{AFKKZ} in order to relate the notion flexibility to the Makar-Limanov and Derksen invariants of an affine variety.

\section{Makar-Limanov invariants}

Given an affine variety $X$, a derivation $\delta$ of $\C [X]$ 
is said to be locally nilpotent if for all $f\in\C [X]$ there exists a non-negative integer $k$ such that $\delta^k(f)=0$. A vector field on an affine variety $X$ can be seen as a derivation $\delta$ on $\C [X]$, and it is well known that its flow defines an  algebraic $\C_+$-action on $X$ if and only if $\delta$ is  locally nilpotent \cite{Fre}.
A regular function $f\in\C [X]$ is $\C_+$-invariant if and only if $\delta (f)=0$.

\begin{definition}
The Makar-Limanov invariant is followig subring of $\C [X]$:

\begin{align}
\rm{ML}(X)=\bigcap_{\delta\in\rm{LND} (X)}\ker \delta 
\end{align}

where $\rm{LND}(X)$ represents the set of all locally nilpotent derivations of $X$.

The field Makar-Limanov invariant is following subfield of $\C (X)$:
\begin{align}
\rm{FML}(X)=\bigcap_{\delta\in\rm{LND} (X)}\ker \delta 
\end{align}
 
where now each $\delta\in \rm{LND} (X)$ is extended as a derivation on the field $\C (X)$ of rational functions of $X$.
We say that the (field) Makar-Limanov invariant is trivial if ($\rm{FML}(X)$) $\rm{ML}(X)=\C$.

\end{definition}

A regular (rational) function belongs to $\rm{ML}(X)$ ($\rm{FML}(X)$) 
if and only if it is invariant 
with respect to all algebraic $\C_+$-actions on $X$. 
Clearly $\rm{ML} (\C^n)\cong \C$, since each partial derivative generates a $\C_+$-action on $\C^n$.
 Makar-Limanov \cite{M-L} introduced this ring in order  to distinguish the Russell threefold
$X=\{ x+x^2y+z^2+t^3=0\}$ from $\C^3$ by proving that $\rm{ML}(X)=\C [x]_{\mid_X}$. In geometrical terms, it means that the orbits of $\C_+$-actions on
the Russell cubic are all tangent to the two-dimensional fibers of the fibration given by the function $x$.

It is clear that if an affine variety has a flexible point, then $\rm{FML}(X)$, and therefore $\rm{ML}(X)$, consists of constant only.
In \cite{AFKKZ} it is proven that  the field $\rm{FML}(X)\cong\C$ if and only if the generic point of $X$ is  flexible, but  the statement at the level of regular functions does not seem to hold:
Liendo \cite{L} constructs an example of a singular affine variety with $\rm{ML}(X)\cong\C$ but with no flexible points. Our next theorem gives a slight improvement on this direction.

\begin{theorem}\label{MLflexible}
Let $X$ be an affine variety. 

\begin{enumerate}

\item $\rm{FML}(X)\cong \C$ if and only if $X$ possess a flexible point if and only if $\saut (X)$ acts with an open orbit;

\item let $X$ have finite Picard group: then $\rm{ML} (X)\cong \C$ if and only if there exists a flexible point of $X$.

\end{enumerate}

\end{theorem}
The first part of the theorem is contained in \cite{AFKKZ}, and we add it here for completeness of explanation.
In order to prove the second part,  we use the following simple lemma and the generalization of Rosenlicht  Theorem on rational invariants stated in Theorem \ref{ratinv}.

\begin{lemma}\label{inv}
Let $X$ be an affine variety. Let $D$ be a hypersurface of $X$, given set theoretically by the zero locus of some function $g\in\C [X]$, and $\delta$ a locally nilpotent derivation on $X$. If $\delta$ is tangent to $D$ for all points of $D$, then 
$\delta (g)=0$.

\begin{proof}
The complement of $D$ is invariant with respect to the $\C_+$-action defined by $\delta$. Hence $g$ restricted to any orbit contained in $X\setminus D$ is a polynomial that does not vanish and therefore constant.

\end{proof}

\end{lemma}

\begin{theorem}\label{ratinv} (\cite{AFKKZ}, Theorem 1.12)
There exists a finite collection of rational $\saut (X)$-invariants which separate 
$\saut (X)$-orbits in general position.
\end{theorem}

For our purpose we need the following corollary of the above theorem.

\begin{corollary}\label{quotient}
There exists an open affine subset $X'\subset X$ which admits a  geometric
quotient $\pi : X'\rightarrow X'//\saut (X)$. 

\begin{proof}
Consider a collection $f_1,f_2,\dots f_N$ of rational $\saut (X)$-invariants as in
Theorem \ref{ratinv}. Let $X'$ be the complement of the polar divisors of 
the functions $f_k$. The algebra $R$ generated by $f_1,\dots , f_N$ is a finitely generated subalgebra of $\C [X']$, and the inclusion of $R$ into $\C [X']$
induces a quotient $X'\rightarrow \Spec R :=X'//\saut (X)$.

\end{proof}

\end{corollary}

\subsubsection*{{\bf Proof of Theorem \ref{MLflexible}}}

\hspace{1cm}

\begin{enumerate}

\item  See \cite{AFKKZ}, Proposition 5.1.

\item Let $\rm{ML} (X)\cong \C$. If we assume that no points are flexible,
then the orbits of $\saut (X)$ have codimension at least one, by part (1) of the Theorem. It follows that the quotient space $X'//\aut (X)$ constructed in Lemma \ref{quotient}
is of dimension at least one. For $x\in X'$,  consider a divisor $\Gamma$ containing the point $\pi (x)$:
then  $D=\pi ^{-1}(\Gamma)$ is a $G$-invariant divisor containing $x$, and it is enough to show the result in the case when the orbit is a codimension one subvariety $D$ containing $x$. In this case, since $X$ has finite Picard group, $D$ as a set is given by the zero locus of some non-costant regular function $g\in \C[X]$, and all derivations are tangent to $D$. Hence, by Lemma \ref{inv}, $\delta (g)=0$ for all locally nilpotent derivations $\delta$, contradicting the triviality of $\rm{ML} (X)$.

\end{enumerate}

\hspace{14.7cm}$\Box $

\hspace{2cm}

\section{Derksen invariant}

We now describe a type of rings which can be naively considered as the `dual" version of the Makar-Limanov ones.  

\begin{definition}
The Derksen invariant of an affine variety $X$ is the subalgebra $\rm{D}(X)$
generated by the collection $\{ \ker \delta \}$, where $\delta$ runs over the set of  non-trivial locally nilpotent derivations on $X$.
Similarly, the field Derksen invariant $\rm{FD}(X)$ is the subfield of $\C (X)$ generated by $\{ \ker\delta\subset \C (X) \}$.
We say that the (field) Derksen invariant is trivial if ($\rm{FD}(X)$) $\rm{D}(X)=\C [X]$.
\end{definition}

Therefore a regular (rational) function belongs to $\rm{D} (X)$ 
($\rm{FD}(X)$) if it is a linear combination of products or regular (rational)
functions that are invariant with respect to at least one non-trivial $\C_+$-action on $X$. It is clear that $\rm{D}(\C^n)$ coincides with the whole polynomial ring in $n$ variables.
In his Ph.D. thesis \cite{D}, Derksen finds a new proof of the fact that the Russel cubic $R$ is not isomorphic to $\C^3$ by showing that 
$\rm{D}(X)\ncong \C [R]$.

\begin{remark}\label{Dnott}
If $X$ is an affine variety of dimension at least two,
and there exists a non-trivial action of $\C_+$ on $X$, then $\rm{D}(X)\ncong\C$.
\end{remark}

It is natural to understand if there is any relation between $\rm{ML}(X)$ and
$\rm{D}(X)$. In particular, one expects that ``small'' $\rm{ML}(X)$ gives us a ``big" $\rm{D}(X)$: however, we have the following examples.

\begin{example}\label{cm}
Crachiola and Maubach \cite{CM} constructed:
\begin{enumerate}

\item
an affine variety $X$ with $\rm{ML}(X)\cong\C$ but $\rm{D}(X)\neq\C [X]$, 
which is not  factorial and not  normal.
\item
an affine variety $X$ with $\rm{D}(X)=\C [X]$ but $\rm{ML}(X)\ncong \C$: 
take $X=\C^*\times \C^2$; $X$ is smooth factorial and $X$ has no flexible points.

\end{enumerate}

\end{example}

\begin{theorem}\label{fld}
Let $X$ be an affine variety with a flexible point.
 Then $\rm{FD}(X)=\C (X)$.
\end{theorem}

The proof is based on the following standard fact of algebraic or differential geometry.

\begin{proposition}\label{birational}

Suppose that there is a finite collection 
$K=\{f_1,\dots f_N\}\subset\C [X]$ such that:

\begin{enumerate}

\item $f_k\in\ker\delta_k$ for some non-trivial locally nilpotent derivation $\delta_k$ on $X$;

\item elements of $K$ separate points on a Zariski open subset $U$ of $X$;

\item for general $x\in X$, the collection $\{ df_k (x) \}$ separates tangent vectors.

\end{enumerate}

Then the map $F:U\rightarrow \C^N$ given by  $F(x)=(f_1(x),\dots ,f_N(x))$
is an open immersion of $U$ into  $\C^N$,  and $\rm{FD}(X)=\C (X)$. If moreover (1)-(3) hold for all points of $X$ and $F$ is a closed embedding of $X$ into $\C^N$, then $\rm{D} (X)=\C [X]$.

%Therefore the functions $f_1,f_2\dots f_N$ generates the algebra of regular functions of the closure of $F(U)$ in $\C^N$, and the field of fractions generated by 
%$f_1,f_2,\dots f_N$ coincides with the field $\C (U)\cong \C (X)$.

\end{proposition}

\begin{convention}
In the rest of the section, we will assume that $X$ is a smooth flexible affine variety: since flexibility is an open condition on the points of $X$, in the proofs we can always replace
$X$ with the flexible affine open subset $U$. 
\end{convention}

\begin{lemma}\label{step1}
For any two distinct points $x$,$y$ in $X$, there is a non-trivial locally nilpotent derivation $\delta$, and regular function $f\in\ker\delta$
that takes different values at $x$ and $y$.
\begin{proof}

We use the $2$-transitivity of the action of $\saut (X)$ on $X$ (Theorem \ref{ft}). 
Let $f$ be any non-costant function in the kernel of some non-trivial locally nilpotent derivation $\delta$.
If $f(x)\neq f(y)$ we are done. If $f(x)=f(y)$, let $z\in X$
be such that $f(x)\neq f(z)$, and consider an automorphism of $X$ such that 
$\psi (y)=y$ and $\psi (x)=z$. Then $f\circ\psi$ belongs to the kernel of 
$\psi_*\delta$ and  takes different values at $x$ and $y$.

\end{proof}
\end{lemma}

\begin{lemma}\label{step3}
There exists a finite collection $K$ satysfing (1) and (2) of Proposition 
\ref{birational}.

\begin{proof}

Let $S=X\times X-\{(x,x)\quad |\quad x\in X\}$.
We denote by $\delta_k$ any non-trivial locally nilpotent derivation,
for $k$ any integer.

Let $(x_1,y_1)\in S $.
Lemma \ref{step1} implies that there is a function $f_1\in\ker\delta_1$ such
that $f_1(x_1)-f_1(y_1)\neq 0$: this inequality holds on a Zariski open 
subset $V_1$ of $S$ containing $(x_1,y_1)$. Let
$(x_2,y_2)$ be a point on an irreducible component of $Z_1=V-V_1$: we can find a function $f_2\in\ker\delta_2$ such that $f_2(x_2)-f_2(y_2)\neq 0$, and this inequality holds on a Zariski open set $V_2$ of $S$ containing $(x_2,y_2)$.  If we repeat this construction, after finitely many steps we produce a Zariski-open covering $\{V_i\}$ of $S$
and regular functions $f_i\in\ker\delta_i$ such that $f_i(x)-f_i(y)\neq 0$ 
for all $(x,y)\in V_i$. Since each $(x,y)\in V_i$ for some $i$, the result follows.

\end{proof}

\end{lemma}

In order to find a collection of invariant functions that separate tangent vectors, we use the following result, which is taken from a generalization of Kleinman's transversality Theorem.

\begin{theorem}\label{kt}(\cite{AFKKZ} Theorem 1.15)
Suppose that $\saut (X)$ acts on $X$ with an open orbit $O$.
Then there are subgroups $H_1,H_2,\dots H_n$ isomorphic to $\C_+$ 
such that for any locally closed 
reduced subscheme $Y$ and $Z$ in $O$ we can find a Zarisky dense open subset
$U\subset H_1\times H_2\times\dots \times H_n$ such that for all $(h_1,h_2,\dots h_n)\in U$ the translate $h_1h_2\dots h_nZ_{reg}$ intersects $Y_{reg}$ transversally.

\end{theorem}

\begin{lemma}\label{tran}

If $X$ has a flexible point, there exists a collection $K'$ that satysfies (1) and (3)
of Proposition \ref{birational}.
\begin{proof}

According to Theorem \ref{kt},
there are two non-trivial locally nilpotent derivations $\delta_1$ and $\delta_2$, and hypersurfaces $D_1=\{f_1=0\}$ and $D_2=\{f_2=0\}$ with $f_k\in \ker\delta_k$ that meet transversally. Consider then $D_1\cap D_2$. Then there exist a locally nilpotent derivation $\delta_3$ and  
a hypersurface $D_3=\{f_3=0\}$, with $f_3\in\delta_3$, that meet $D_1\cap D_2$ transversally. 
In particular $df_1\wedge df_2\neq 0$ for a general point of the intersection of the two subvarieties, hence on a Zariski open subset of $X$.
Let $x$ be a point where $D_1\cap D_2$ intersect $D_3$ transversally. 
Then there exists a $v\in T_x(D_1\cap D_2) $ such that 
$df_3 (v)\neq 0$, since the tangent space $T_x(D_1\cap D_2)$
is not contained in $T_xD_3$. Hence $df_3(x) $ is not a linear combination of 
$df_1(x)$ and $df_2(x)$ (such a vector $v$ is annihilated by the differentials
$df_1(x)$ and $df_2(x)$, if $df_3(x)$ were such a combination then $df_3(x)(v)=0$). Apply Kleinman's theorem inductively, we find a collections of non-trivial locally nilpotent derivations $\delta_k$ and hypersurfaces 
$D_i=\{f_i=0\}$, $i\in\{1,\dots n\}$, such that 
$df_1\wedge df_2\wedge\dots\wedge df_n\neq 0$ on a Zariski open set of $X$.

\end{proof}

\end{lemma}

\subsubsection*{{\bf Proof of Theorem \ref{fld}}}

\hspace{1cm}

The collection $K\cup K'$, where
$K$ is from Lemma \ref{step3} and $K'$ is from Lemma \ref{tran}, satisfyies (1)-(3) of Proposition \ref{birational}.

\hspace{14.7cm} $\Box$

Theorem \ref{MLflexible} and Theorem \ref{fld} togheter imply the following.

\begin{corollary}
Let $X$ be an affine variety with finite Picard group. 
Then $\rm{ML}(X)\cong\C$ implies $\rm{FD}(X)\cong\C (X)$

\end{corollary}

%\begin{remark}
%Except for the two dimensional case, that will be discussed in the next section, we do not know how to construct a closed embedding satisfying the conditions of Proposition \ref{birational},
%and therefore replace $\rm{FD}(X)$ with $\rm{D}(X)$ in the above corollary.
%\end{remark}

\section{Affine surfaces}

In this section we specialize our previous results to the case of surfaces, in which we can obtain stronger conclusions.

\begin{theorem}\label{flexsur}
If $X$ is a smooth flexible affine surface. Then 
then $\rm{D}(X)=\C [X]$  and $\rm{ML}(X)\cong\C$.

\begin{proof}
It is clear that $\rm{ML} (X)\cong \C$.

As for the Derksen invariant, consider the embedding $F:X\rightarrow \C^N$ constructed in Proposition
\ref{birational}: that is $F=(f_1,f_2,\dots f_N)$ for $f_k\in\Ker \delta_k$, where $\delta_k$ is some non-trivial locally nilpotent derivation.
We want to show that the image of $X$ under $F$ is closed, in order to apply the last sentence of Proposition \ref{birational}.
Consider the Zariski closure $X'$ of the image $F(X)$ in $\C^N$. 
We show that if $X'-F(X)$ is not empty, then we can construct a compact curve inside the affine variety $X'$, which is impossible.
Any irreducible component of $X'-F(X)$ is a curve: if not, we can apply Hartog's Theorem to the normalization of $X'$ to run into a contradiction.
%Suppose first that $D$ is a point of coordinates $(x^0_1,x^0_2,\dots ,x^0_n)$.
%The extension of $f_k$ to $X'$ coincides with the coordinate function $x_k$ given by the inclusion in $\C^N$. Therefore the equation $f_k=x^0_k$ gives a 
%divisor in $X'$ which by Lemma \ref{lines} consists of a line in $X$ plus the extra point $D$,
%that is a compact curve.
%Suppose now that $D$ is a curve. Then at least one coordinate function, let's say $x_1$, is not constant. Therefore for a generic value $x_1^0$
%the divisor $f_1=x^0_1$ is densely contained in $X$, and its closure is again a compact curve in the affine variety $X'$.
Consider then a general point $p$ that belongs to a unique irreducible component $D$  of the boundary, with coordinates $(c_1,c_2,\dots c_n)$.
Since along $D$ at least one coordinate, say $x_1$, is non-constant the curve $\Gamma$ of equation $f_1=c_1$ must intersect $D$ at $p$ transversally.
If we apply the next lemma, we conclude that $\Gamma \cap X$ is a line: this is impossible since  $\Gamma$ would be compact inside $X'$.

\end{proof}

\end{theorem}

\begin{lemma}\label{lines}
Let $\delta$ be a non-trivial locally nilpotent derivation on $X$. 
Let $f\in\Ker\delta$. then the levels sets $f=c$ consists of union of lines. 

\begin{proof}
The results is a special case of (\cite{GMMR} Lemma 1.6), where there is no need of assumptions on the Picard group, but we provide a proof here.

The level curves of the $\C_+$-action associated to $\delta$ are union of orbits
. If $x$ is not a fixed point, then $\C.x$ is a line, therefore we are left to understand the fixed points.
Since the fixed points of a locally nilpotent derivation can not be isolated 
(see \cite{Fre}, Theorem 1.30), $F$ is a curve. Since $X$ has finite Picard group, $F$ is set-theoretically the zero locus of a regular function $f$ and 
and by  Lemma \ref{inv} we obtain that $\delta (f)=0$.
Moreover, we can write $\delta=f^k\Delta$, where $\Delta $ is a non-trivial derivation on $F$. For $c\neq 0$, on the smooth locus of the level curves $f=c$,
we have $df(\Delta)=\Delta (f)=0$, and by continuity this condition holds
also on $F$. That is $\Delta $ is tangent to $F$. Since $\Delta$ is locally nilpotent, the restriction to $\Delta$ on $F$ gives a non-trivial $\C_+$-action,
hence each irreducible component of $F$ consists of a line (\cite{Fre}, Corollary 1.29). The components are pairwise disjoint, since each of them is an orbit.

\end{proof}
\end{lemma}

We finish up by applying the above constructions to some examples of affine surfaces with trivial Makar-Limanov invariant, that have been studied by many authors.

\vspace{0.5cm}

{\bf I. $\mathbb{Q}$-homology planes}

A $\mathbb{Q}$-homology plane is a smooth affine surface with the same rational homology of the affine plane.
A smooth affine surface with trivial Makar-Limanov invariant and  finite Picard group is a $\mathbb{Q}$-homology plane (\cite{GMMR}, Lemma 1.4).
Here we will use the following result of Masuda and Miyanishi   to prove that $\rm{D}(X)=\C [X]$ for a $\mathbb{Q}$-homology plane,
 a result that improves Lemma 5.3 of \cite{GMMR}.

\begin{theorem}(\cite{MM}, Theorem 3.4)
Let $X$ be a $\mathbb{Q}$-homology plane with $\rm{ML}(X)\cong \C$. Then $Pic (X)\cong \mathbb{Z}_m$, and $X$ is isomorphic to the quotient
of the affine surface $Y=\{ xy=z^m-1\}$ with respect to a free $\mathbb{Z}_m$-action 
$\tau(x,y,z)=(\tau x, \tau^{-1}y, \tau^jz)$ ($j<m$ and $\gcd (m,j)=1$).
\end{theorem}

\begin{theorem}\label{QHP}
Let $X$ be a  $\mathbb{Q}$-homology plane with $\rm{ML}(X)\cong\C$. Then $\rm{D}(X)\cong \C [X]$.
\begin{proof}
The following locally nilpotent derivations on the surface $Y$ 

\begin{align}
\nn \delta_1= x^{j-1}\left( mz^{m-1}\frac{\partial}{\partial y}+x\frac{\partial}{\partial z}\right)\\
\nn \delta_2=y^{j-1}\left( mz^{m-1}\frac{\partial}{\partial x}+y\frac{\partial}{\partial z}\right)
\end{align}

are are $\mathbb{Z}_m$-invariant, hence they descends to locally nilpotent derivation on $X$.
Their common zero locus consists of one point on $X$, namely the image $p$ of the  set $\{x=y=0\}$  on Y.  Moreover, the derivations are transversal at points other than $p$, or, using the terminology  of Section 1, $X-p$ is flexible. 
The map $F:(X-p)\rightarrow \C^N$ of Proposition \ref{birational} extends to $X$, is injective and separates all tangent vectors in $X-p$. To separate the tangent vectors at $p$, 
observe that the zero loci of the invariant functions $x^m$ and $y^m$ are transversal at $p$, 
therefore the map $(F,x^m,y^m):X\rightarrow \C^{N+2}$ becomes a closed embedding.
\end{proof}
\end{theorem}

\begin{remark}
The author is not aware whether $p$ is a fixed point of the action of the whole group of authomorpshims of $X$. 
If it is not the case, then $X$ is a flexible point and the previous theorem is a consequence of Theorem \ref{flexsur}. 
\end{remark}

\vspace{1cm}

{\bf II. Danilov-Gizatullin surfaces}

Given a Hirzebruch surface $\Sigma_d$, by definition a Danilov-Gizatullin surface $V$ is the complement of an ample section $S$ 
of a ruling of $\Sigma_d$. The Danilov-Gizatullin Theorem states that the isomorphism class of $V=V_n$ depends only on $n=S.S$,
 the self intersection number of $S$ (the ampleness of $S$ imposes that $n>2$ and $n>d$).

In \cite{P}, Example 1.19, Popov gives a list all affine surfaces which admit an action of an algebraic group $G$ with an 
 open orbit whose complement consists of at most finitely many points. Such surfaces are said to be quasi-homogeneous in the sense of Gizatullin.
Among the Danilov-Gizatullin surfaces, only $V_2$ appears in his list.
Hence we can prove the following fact:

\begin{lemma}
Let $n>2$. There is no non-trivial $SL_2$ or $PSL_2$-actions on $V_n$.
\begin{proof}
Let $G\cong SL_2$ act on  $V_n$. For what just observed , if there is an 
$SL_2$-orbit $O$ of maximal dimension, its complement $D$ must be of dimension $1$.
But $D$ must contain a $1$-dimensional orbit $\Gamma$ : from the standard theory of reductive groups, the closure of $O$ contains a unique closed orbits, and fixed points are closed orbits.
Therefore $\Gamma\cong SL_2/I$, and $I$ is a two-dimensional algebraic subgroup of $SL_2$. Now the claim is that the connected component $I_0$ of $I$ is parabolic: therefore $\Gamma\cong (SL_2/I_0)/(I/I_0)$
is a projective variety contained in $V_n$, impossible.

\end{proof}

\end{lemma}

Therefore, the surfaces $V_n$ for $n>2$ do not belong to the class of surfaces satisfying the following Proposition: 

\begin{proposition}(\cite{P}, Theorem 1.25 and Remark 1.26)
If $\aut(X)$ contains an $SL_2$ or $PSL_2$, then $D(X)=\C [X]$
\end{proposition}

However,  the Danilov-Gizatullin surfaces are flexible and  Theorem \ref{flexsur} applies.

\begin{theorem}
For any $n$, the Danilov-Gizatullin surfaces have trivial Makarl-Limanov and Derksen invariants:
that is $\rm{D}(V_n)= \C [V_n]$ and $\rm{ML}(V_n)\cong\C$.

\end{theorem}

In alternative, we can easily calculate directly the Derksen invariant for 
Danilov-Gizatullin surfaces (see the Appendix).

\vspace{1cm}

{\bf III. Quasi-homogeneous surfaces}

In Section 5 of \cite{GMMR}, the authors claim that no examples of surfaces with trivial Makar-Limanov invariant but non-trivial Derksen invariant are known.
Possible candidates are given by the following result.

\begin{theorem}\label{qhom}
Let $X$ be a surface with trivial Makar-Limanov invariant 
such that the complement of the open orbit $O$ of $\aut (X)$ consists of a finite
set of cardinality at least two. Then $\rm{D}(X)\neq\C [X]$.

\begin{proof}
Let $A$ be the complement of the open orbit $O$ of $\aut (X)$.
Since the points in $O$ are all flexible, 
Proposition \ref{birational} gives an embedding 
$F:X-A\rightarrow \C^N$. Moreover, 
$F$ extends as a regular map on $A$ since $X$ is smooth.
Using the same reasoning of Theorem \ref{flexsur}, we conclude that  
$F(X)$ is closed in $\C^N$. Using the notation of Proposition \ref{birational},
consider the first component $f_1$ of $F$. Let $f_1\in\Ker\delta_1$,
and assume that $\delta_1$ is irreducible (that is, we can not write $\delta=g\delta'$ for some $g\in\Ker\delta$). 
The $\C$-fibration $X\rightarrow X//\delta$ has only one singular fiber
(\cite{GMMR}, Lemma 1.6), hence the set $A$ is contained in one level set of $f_1$. The same reasoning applies to all components $f_k$ of $F$:
hence the image of $A$ under $F$ consists of a unique point $P$ in $X'=F(X)$,
Moreover, the map $F:X\rightarrow X'$ is proper and quasifinite 
(that is finite): hence $X'$ is not normal at $P$, otherwise $F$ will be an isomorphism. 
Consider now a regular function $g\in \C [X]$ which is not regular in $X'$. 
We claim that $g\notin \rm{D}(X)$. Suppose that instead $g$ can be written
 as a combination of products of invariant functions $f_{N+1},f_{N+2}, \dots f_{M}$.
Then the map $\tilde {F}=(f_1,f_2,\dots f_M) :X\rightarrow \C^{M}$ is a closed embedding that factors through $F$ 
and its image $\tilde{X }$ is isomorphic to $X'$ via the projection onto the first $N$ factors. Therefore $g \in \C [\tilde{F} (X)]$
descends to a regular function on $X'$, contradiction.

\end{proof}

\end{theorem}

\begin{corollary}
Let $X$ be a surface with trivial Makar-Limanov invariant. If it has trivial Derksen invariant, then the action of $\Aut (X)$ on $X$ has at most one fixed point.
\end{corollary}

\section{Appendix: direct computation of $\rm{D}(V_n)$}
In alternative, we can calculate directly the Derksen invariant for 
Danilov-Gizatullin surfaces, using the following Proposition \cite{Do}:

\begin{proposition}
$\C [V_n]$ is generated by the monomials

\begin{align}
&x_k=a_2^{b-k}a_3a_4^k\\
&y=a_1a_2
\end{align}

\end{proposition}

One can check that the following are well-defined locally nilpotent derivations
on $V_n$:

\begin{align}
&\delta =ba_2^{b-1}a_3\frac{\partial}{\partial a_1}+a_4\frac{\partial}{\partial a_2}\\
&\delta'=a_1^{b-1}a_2^{b}\frac{\partial}{\partial a_4}+a_1^{b}\frac{\partial}{\partial a_3}
\end{align}
where $b=n-1$ and that:

\begin{align}
&\delta (y)=1+nx_0\\
&\delta (x_k)=(b-k)x_{k+1}\\
&\delta'(y)=0
\end{align}

The $\C_+$-action defined by the derivation $\delta$ acts as follows:

\begin{align}
e^{t\delta}(y)=y+t(1+nx_0)+t^2\frac{nb}{2}x_1+\dots +t^{b+1}\frac{nb!}{(b+1)!}x_b
\end{align}

For here we conclude that 

\begin{align}
p_t=tx_0+t^2\frac{b}{2}x_1+\dots + t^{b+1}\frac{nb}{(b+1)!}x_b
\end{align}
belong to $\rm{D}(V_n)$.

For pairwise distinct non zero numbers $t_0,t_2\dots t_b$
the matrix

\begin{align}
\left(
\begin{array}{cccc}
t_0 & t_0^2 & \ldots & t_0^{b+1}\\
t_1 & t_1^2 & \dots & t_1^{b+1}\\
\cdots & \cdots & \cdots & \cdots\\
t_b & t_b^{2} & \ldots & t_b^{b+1}
\end{array}
\right)
\end{align}

has non-zero determinant. 
Therefore $p_{t_k}$ for a basis for the vector space spanned by the monomials
$x_0,x_2\dots x_{b}$, and the result follows.

\subsection*{Acknowledgments}
 I would like to thank: Shulim Kaliman, Dror Varolin and Ivan Arzhantsev for reading preliminary versions of the paper and the useful conversations;
Anthony Crachiola, Steven Lu, Alok Maharana and Stefan Maubach for interesting discussions related to the subjet.

\vfuzz=2pt
%%%%%%%%%%%%%%%%%%%%%%%%%%%%%%%%%%%%%%%%%%%%%%%%%%%%%%%%%%%%%%%%%%%
%%%%%%%%%%%%%%%%%%%%%%%%%%%%%%%%%%%%%%%%%%%%%%%%%%%%%%%%%%%%%%%%%%%
%%%%%%%%%%%%%%%%%%%%%%%%%%%%%%%%%%%%%%%%%%%%%%%%%%%%%%%%%%%%%%%%%%%

\providecommand{\bysame}{\leavevmode\hboxto3em{\hrulefill}\thinspace}

\end{document}